\documentclass[a4paper]{mathscan}

\usepackage{latexsym}
\usepackage{amssymb}
\usepackage[numbers]{natbib}

\newtheorem{theorem}{Theorem}[section]
\newtheorem{lemma}[theorem]{Lemma}
\newtheorem{proposition}[theorem]{Proposition}
\newtheorem{corollary}[theorem]{Corollary}

\theoremstyle{definition}
\newtheorem{definition}[theorem]{Definition}

\theoremstyle{remark}
\newtheorem{remark}[theorem]{Remark}
\newtheorem{example}[theorem]{Example}

\numberwithin{equation}{section}

\newcommand{\C}{\mathbb{C}}
\newcommand{\R}{\mathbb{R}}
\newcommand{\Z}{\mathbb{Z}}
\newcommand{\N}{\mathbb{N}}
\newcommand{\E}{\mathbb{E}}

\newcommand{\re}{\mathrm{Re}}

\newcommand{\f}{\varphi}
\newcommand{\supp}{\mathrm{supp}}

\newcommand{\om}{\omega}
\newcommand{\Hh}{\mathcal H}

\begin{document}

\title[Relative inner amenability]{Relative inner amenability and relative property gamma}

\author{Paul Jolissaint}
\address{Universit\'e de Neuch\^atel, Institut de Math\'ematiques, Rue Emile-Argand 11, CH-2000 Neuch\^atel, Switzerland}
\email{paul.jolissaint@unine.ch}


\date{\today}


\begin{abstract}
Let $H$ be a proper subgroup of a discrete group $G$. We introduce a notion of relative inner amenability of $H$ in $G$, we prove some equivalent conditions and provide examples coming mainly from semidirect products, as well as counter-examples. We also discuss the corresponding relative property gamma for pairs of type II$_1$ factors $N\subset M$ and we deduce from this a characterization of discrete, icc groups which do not have property (T).
\end{abstract}

\maketitle

\section{Introduction}\label{intro}

The aim of the present notes is to introduce a relative version of inner amenability for pairs of groups $H\subset G$ and a relative version of property gamma for pairs of type II$_1$ factors $N\subset M$ with separable preduals.

Inner amenability was first introduced by E. Effros in \cite{Effros} and was further studied by W. Paschke \cite{Paschke}, then by E. B\'edos and P. de la Harpe in \cite{BH}, and in \cite{HS} by P. de la Harpe and G. Skandalis. The aim of E. Effros was to translate Murray and von Neumann's \textit{property gamma} of type II$_1$ factors into a property for groups. More precisely, assume that $H$ is an icc group (\textit{i.e.} every non trivial conjugacy class is infinite) and that its group von Neumann algebra $L(H)$ has property gamma (see definition below). Then $H$ is 
\textit{inner amenable} in the sense that the C$^*$-algebra $\ell^\infty(H\smallsetminus\{1\})$ has an inner invariant state (equivalently, the $H$-set $H\smallsetminus \{1\}$ has an invariant mean). The converse remained open for almost forty years and was proved to be false by S. Vaes in \cite{Vaes}.

Here we consider the situation where $H$ is a proper subgroup of a group $G$. As the subset $G\smallsetminus H=\{g\in G : g\notin H\}$ is invariant under conjugation by elements of $H$, \textit{i.e.} $g(G\smallsetminus H)g^{-1}=G\smallsetminus H$ for every $g\in H$, this leads to the following natural definition.

\begin{definition}\label{def1.1}
(1) Let $H\subset G$ be a pair of groups. Then we say that $H$ is \textbf{inner amenable relative to} $G$ if $H$ is a proper subgroup of $G$ and if
there exists an $H$-invariant state on $\ell^\infty(G\smallsetminus H)$ for the action of $H$ on $G\smallsetminus H$ by conjugation.\\
(2) Let $N\subset M$ be type II$_1$ factors. Then we say that $N$ has \textbf{property gamma relative to} $M$ if, for every finite set $F\subset N$, there exists a bounded sequence $(x_n)\subset M$ such that $\Vert x_n\Vert_2=1$ and $\E_N(x_n)=0$ for every $n$, and $\Vert yx_n-x_ny\Vert_2\to 0$ for every $y\in F$.
\end{definition}

In Section \ref{InnerAm}, we give a few equivalent conditions of relative inner amenability; most of them are reformulations of amenability of actions of groups on sets. Here is a sample of conditions equivalent to relative inner amenability:

\par\vspace{2mm}\noindent
\textbf{Theorem A.}
\textit{Let $H$ be a proper subgroup of the group $G$. Then the following conditions are equivalent:}
\begin{enumerate}
\item [(1)] \textit{$H$ is inner amenable relative to $G$.}
\item [(2)] \textit{There exists a net $(\xi_n)$ of unit vectors in $\ell^2(G)$ such that $\supp(\xi_n)\subset G\smallsetminus H$ for every $n$ and
\[
\lim_{n}\Vert\alpha(h)\xi_n-\xi_n\Vert_2=0
\]
for every $h\in H$, where $\alpha$ denotes the representation by conjugation of $G$.}
\item [(3)] \textit{There exists a state $\f$ on $B(\ell^2(G))$ such that $\f(\alpha(h))=1$ for every $h\in H$ and $\f(e_H)=0$, where $e_H$ is the projection onto $\ell^2(H)$.}
\end{enumerate}


\begin{remark}\label{rem2.3}
(1) A trivial situation where a proper subgroup $H$ of a group $G$ is inner amenable is when $H$ is an amenable group: it is inner amenable relative to any group $G$ containing it. Indeed, for instance fix any $g_0\in G\smallsetminus H$ and let $X=\{hg_0h^{-1} : h\in H\}$. Then $h\mapsto hg_0h^{-1}$ is an $H$-map from the amenable $H$-space $H$ to $X$, and Lemma \ref{lem1.4} applies.

Another rather trivial situation is when 
there exists some element $g\in G\smallsetminus H$ such that the corresponding orbit $\{hgh^{-1} : h\in H\}$ is finite. If it is the case, we will say that $H$ is \textit{trivially inner amenable relative to} $G$. It holds for instance when the group $G$ is a direct product $G=H\times K$ with any non-trivial group $K$.

Therefore, in order to discuss interesting instances, we introduce the following condition:
\begin{center}
$(\star)$ $\quad$ $\{hgh^{-1}:h\in H\}$ \textit{is infinite for every} $g\notin H$.
\end{center}
(2) Observe that if $H$ has Kazhdan's property (T) and if it is a proper subgroup of $G$ then it is inner amenable relative to $G$ if and only if it is trivially inner amenable relative to $G$. This follows from \citep[Lemma 4.2]{GlasMon}. 

More generally, let $K\subset H\subset G$ be three groups so that $H$ is inner amenable relative to $G$ and that the pair $K\subset H$ has property (T) (\cite{jolTpairs}). By Example \ref{ex2.6}, $K$ is inner amenable relative to $G$, but property (T) of the pair $K\subset H$ implies the existence of some $g\in G\smallsetminus H\subset G\smallsetminus K$ such that $\{kgk^{-1} : k\in K\}$ is finite. In other words, even if $H$ is non-trivially inner amenable relative to $G$, $K$ is trivially inner amenable relative to $G$.
\end{remark}

These observations lead to the following natural question:

\par\vspace{2mm}\noindent
\textbf{Question.} \textit{Which groups $H$ can be embedded into groups $G$ in such a way that they satisfy condition $(\star)$ and
such that $H$ is inner amenable relative to $G$?}

\par\vspace{2mm}
Partial answers are given in Section \ref{InnerAm}. Observe that by Remark \ref{rem2.3}, a necessary condition on $H$ is that it does not have property (T), but we will see in Example \ref{ex2.12} that there are non-Kazhdan groups for which there exists no group $G$ satisfying the conditions of the above question.

The following theorem shows that infinite groups acting amenably on infinite sets with only infinite orbits provide a family of groups that satisfy conditions of the above question.

\par\vspace{2mm}\noindent
\textbf{Theorem B.} 
\textit{Let $H\curvearrowright X$ be an amenable action in the sense of Section \ref{prerequisites} and let $Z$ be any non-trivial group. Let $G=Z\wr_XH$ be the corresponding restricted wreath product group. If all orbits of $H\curvearrowright X$ are infinite, then $H$ is non-trivially inner amenable relative to $G$.}

\begin{remark}\label{rem2.9}
The article of Y. Glasner and N. Monod \cite{GlasMon} provides numerous examples of countable, non-amenable groups satisfying all conditions of Theorem B. 

Following \cite[Definition 1.3]{GlasMon}, let us denote by $\mathcal A$ the class of all countable groups $H$ that admit a faithful, transitive, amenable action on some countable set, and let us say that a countable group $H$ has \textit{property} (F) if any amenable $H$-action (on a countable set) has a fixed point.

Then, by \citep[Theorem 1.5]{GlasMon}, if $H$ and $K$ are countable groups, their free product $H*K$ belongs to $\mathcal A$ unless $H$ has property (F) and $K$ has virtually property (F) (possibly upon exchanging $H$ and $K$). For example, $H*K$ is in $\mathcal A$ as soon as one of the groups is residually finite or non-finitely-generated or amenable.

The authors of \cite{GlasMon} also introduce the following class of infinite, countable groups (\citep[Definition 4.1]{GlasMon}): denote by $\mathcal B$ the class of all countable groups admitting some amenable action on a countable set without finite orbits. Obviously, $\mathcal B$ contains $\mathcal A$, but it is much wider: for instance, any group with a quotient in $\mathcal B$ belongs to $\mathcal B$. Moreover, \citep[Lemma 2.16]{GlasMon} shows that any non-finitely generated, countable group belongs to $\mathcal B$.
\end{remark}

\begin{remark}
It would be interesting to know whether $\mathcal B$ contains all infinite, countable groups that have some weak form of amenability, for instance the Haagerup property \cite{ccjjv}. 

Other weak forms of amenability are, on the one hand, \textit{weak amenability} due to Cowling and Haagerup \cite{CH}, and, on the other hand, the more recent \textit{weak Haagerup property} introduced by S. Knudby in \cite{Knudby} which generalizes the Haagerup property and weak amenability. As these classes contain (infinite) groups with property (T), they cannot be contained in $\mathcal B$.
\end{remark}

Section \ref{InnerAm} ends with a discussion of pairs $H<G$ where $H$ has finite index in $G$: we show that there are pairs such that $H$ is not inner amenable relative to $G$, even if $H$ itself is inner amenable. This has to be compared to \citep[Ajout]{BH} where the authors prove that if $H$ is inner amenable and if it is of finite index in $G$, then $G$ is inner amenable.

\par\vspace{2mm}
Section \ref{propgamma} is devoted to a study of relative property gamma for pairs of type II$_1$ factors $N\subset M$ with separable preduals. Notice that this relative property is different from Bisch's relative property gamma for the inclusion $N\subset M$ as studied in \cite{Bi1} and \cite{Bi2}: indeed, D. Bisch says that the \textit{inclusion} $N\subset M$ has \textit{property gamma} if there is a sequence of unitary elements $(u_n)\subset U(N)$ such that $\tau(u_n)=0$ for every $n$ and which is central for $M$, \textit{i.e.} if $\Vert xu_n-u_nx\Vert_2\to 0$ for every $x\in M$. Murray and von Neumann's property gamma corresponds to the case $N=M$.

Because of the separability condition on the preduals of $N$ and $M$, both properties can be expressed in terms of relative commutants in ultraproduct algebras (see below for precise definitions): if $N$ is a subfactor of the type II$_1$ factor $M$, one has the following chain of natural inclusions:
\[
M'\cap N^\om\subset N'\cap N^\om\subset N'\cap M^\om.
\]
Then our relative property gamma means that $N'\cap M^\om$ contains strictly the relative commutant $N'\cap N^\om$, and Bisch's property gamma for the inclusion $N\subset M$ is equivalent to the non-triviality of $M'\cap N^\om$. We observe that there is no obvious relationship between D. Bisch's notion and ours.

We say that $N$ is \textit{irreducible in} $M$ if $N'\cap M\subset N$, which is equivalent to the equality
$N'\cap M=Z(N)$ where the latter denotes the center of $N$. 

\par\vspace{2mm}

The next theorem, contained in Section \ref{propgamma}, is a characterization of icc groups which do not have Kazhdan's property (T) in terms of relative property gamma for their group von Neumann algebra.

\par\vspace{2mm}\noindent
\textbf{Theorem C.}
\textit{(1) Let $N$ be a type II$_1$ factor with property (T) in the sense of A. Connes and V. Jones \cite{CJ}. If $N$ is an irreducible subfactor of some type II$_1$ factor $M$, then $N$ does not have property gamma relative to $M$.\\
(2) Let $G$ be a countable icc group which does not have property (T). Then its group von Neumann algebra $L(G)$ can be embedded into a type II$_1$ factor $M$ as an irreducible subfactor so that $L(G)$ has property gamma relative to $M$. \\
More precisely, if $G$ does not have property (T), then 
there exists an action $\sigma$ of $G$ on the hyperfinite II$_1$ factor $R$ such that the von Neumann algebra $L(G)$ has the following properties:}
\begin{enumerate}
\item [(a)] \textit{$L(G)$ is an irreducible subfactor of the crossed product $R\rtimes_\sigma G$;}
\item [(b)] \textit{$L(G)$ has property gamma relative to $R\rtimes_\sigma G$.}
\end{enumerate}

\begin{remark}
In some sense, Theorem C shows the following fact: the countable, icc, non-Kazhdan's groups are precisely those groups whose von Neumann algebra embeds irreducibly into type II$_1$ factors and possesses relative property gamma.
\end{remark}

Section \ref{propgamma} also contains an example borrowed from \cite{Vaes} mentioned above: let $H$ be an icc group contained in a countable group $G$ so that $L(H)$ is irreducible in, and has property gamma relative to $L(G)$. Then it is obvious that $H$ is inner amenable relative to $G$, but the converse is false. We use S. Vaes' example to construct a pair $H\subset G$ such that $H$ is inner amenable relative to $G$, but $L(H)$ does not have property gamma relative to $L(G)$.

\par\vspace{2mm}
Section \ref{propgamma} ends with a discussion of pairs $N\subset M$ where $N$ is an irreducible subfactor of $M$ with finite index.

\par\vspace{2mm}
\textit{Acknowledgements.} I am very grateful to Pierre de la Harpe for his very careful reading of a previous version of my article and for his enlightening comments. I also thank Alain Valette for his contribution in Example \ref{ex2.5} and the referee for his careful reading of the manuscript and his valuable comments.

\section{Prerequisites on amenable actions and finite von Neumann algebras}\label{prerequisites}

The present section is devoted to fixing our notation and reminding some facts on:
\begin{itemize}
\item amenable actions of groups on sets, in the sense of the existence of invariant means (or states on $\ell^\infty$);
\item finite von Neumann algebras and their ultraproducts.
\end{itemize}

\par\vspace{2mm}

Let $G$ be a group acting on a set $X$; we denote it by $G\curvearrowright X$.
The corresponding unitary representation $\pi_X:G\rightarrow U(\ell^2(X))$ is defined by
\[
(\pi_X(g)\xi)(x)=\xi(g^{-1}x)
\]
for every $\xi\in \ell^2(G)$ and all $g\in G$, $x\in X$. Observe that the dense subspace
$\ell^1(X)$ is invariant under $\pi_X$. We denote by $\ell^1(X)^+$ the cone of elements $\eta\in\ell^1(X)$ such that $\eta(x)\geq 0$ for all $x\in X$.

In the particular case where $X=G$, one considers the following three actions $G\curvearrowright G$: the action by left translation $g\cdot g'=gg'$, the action by right translation
$g\cdot g'=g'g^{-1}$ and the action by conjugation $g\cdot g'= gg'g^{-1}$. Their associated representations on $\ell^2(G)$ are denoted by $\lambda,\rho$ and $\alpha$ respectively. Explicitly,
\[
\begin{array}{lll}
(\lambda(g)\xi)(g') &=& \xi(g^{-1}g')\\
(\rho(g)\xi)(g') &=& \xi(g'g)\\
(\alpha(g)\xi)(g') &=& \xi(g^{-1}g'g)
\end{array}
\]
for all $g,g'\in G$ and all $\xi\in\ell^2(G)$.
\par\vspace{2mm}
Let $S$ be a subset of a set $X$; we denote by $\chi_S$ its characteristic function, and, when $S$ is finite,
by $|S|$ its cardinal. For any function $f$ defined on some set $Y$ and which has all its values in some group $\mathcal G$ with neutral element denoted by $e$, we denote by $\supp(f)$ the set of elements $y\in Y$ such that $f(y)\not=e$. Finally, for a subset $E$ of $\mathcal G$, we write $E^*:=E\smallsetminus\{e\}$.

\par\vspace{2mm}
Recall that the trivial representation of $G$ is \textit{weakly contained} in a unitary representation $(\pi,\Hh)$ if there exists a net of unit vectors $(\xi_n)\subset \Hh$ such that 
\[
\lim_{n}\Vert \pi(g)\xi_n-\xi_n\Vert=0
\]
for every $g\in G$. Notice that it is equivalent to say that the net of positive definite functions $(\f_n)$, defined by $\f_n(g)=\langle \pi(g)\xi_n|\xi_n\rangle$ for every $g$, converges pointwise to the constant function $1$.
\par\vspace{2mm}

Let $(\pi,\Hh)$ be a unitary representation of $G$; for $f\in\ell^1(G)$, we denote by $\pi(f)$ the associated operator acting on $\Hh$:
\[
\pi(f)=\sum_{g\in G}f(g)\pi(g).
\]

\par\vspace{2mm}
The following lemmas are essentially well-known and generalize the classical case of amenable groups. See for instance \citep[Theorem 4.1]{Gre},
\citep[Theorem 1.1]{KechTsan} or \cite{Rose}. 

We omit almost all proofs, except for the last condition in Lemma \ref{lem1.2}, which is largely inspired from \citep{dlHRV}, and the one-line proof of Lemma \ref{lem1.4}.

\begin{lemma}\label{lem1.2}
For an action $G\curvearrowright X$, the following conditions are equivalent:
\begin{enumerate}
\item [(1)] There exists a $G$-invariant state on $\ell^\infty(X)$, i.e. there exists a positive linear functional $\f$ on $\ell^\infty(X)$ such that $\f(1)=1$ and $\f(g\cdot a)=\f(a)$ for every $a\in\ell^\infty(X)$ and every $g\in G$. Equivalently, there exists a $G$-invariant mean $\mu$ on the set of all subsets of $X$.
In other words, the action $G\curvearrowright X$ is \textbf{amenable} in the sense of \citep[Definition 1.1]{GlasMon}.
\item [(2)] There exists a net of unit vectors $(\eta_n)\subset \ell^1(X)^+$ such that
\[
\lim_{n}\Vert \pi_X(g)\eta_n-\eta_n\Vert_1=0
\]
for every $g\in G$.
\item [(3)] There exists a net of unit vectors $(\xi_n)\subset\ell^2(X)$ such that 
\[
\lim_{n}\Vert \pi_X(g)\xi_n-\xi_n\Vert_2=0
\]
for every $g\in G$; in other words, the trivial representation $1_G$ is weakly contained in $\pi_X$.
\item [(4)] There exists a net $(F_n)$ of non-empty finite sets such that $F_n\subset X$ for every $n$ and
\[
\lim_{n}\frac{|F_n\bigtriangleup g F_n|}{|F_n|}=0
\]
for every $g\in G$. These sets are called \textbf{F\o lner sets}.
\item [(5)] For every finite set $F\subset G$ and for every function $f:F\rightarrow [0,1]$ such that $\sum_{g\in F}f(g)=1$, the number $1$ belongs to the spectrum of $\pi_X(f)$.
\end{enumerate}
\end{lemma}
\textsc{Proof.} Remark that (3) $\Rightarrow$ (5) is easy, so we only need to prove that (5) $\Rightarrow$ (3): It is inspired by \citep[Lemma 3]{dlHRV}. Observe first that (3) is equivalent to the fact that, for every finite, symmetric set $F\subset G$, one has $\kappa(F)=0$, where
\[
\kappa(F):=\inf_{\xi,\Vert\xi\Vert_2=1}\max_{g\in F}\Vert \pi_X(g)\xi-\xi\Vert_2.
\]
Hence, if condition (3) does not hold, there is a finite, symmetric set $F\subset G$ such that $\kappa(F)>0$. Let us prove that $1$ does not belong to the spectrum of the selfadjoint operator $a=\frac{1}{|F|}\sum_{g\in F}\pi_X(g)$. Set $2\delta:=\kappa(F)>0$ and $n=|F|$. Observe that $\Vert a\Vert\leq 1$. 
Then, for every unit vector $\xi\in\ell^2(X)$, there exists $g\in F$ such that $\Vert\pi_X(g)\xi-\xi\Vert\geq\delta$. This implies that $\re\langle \pi_X(g)\xi|\xi\rangle\leq 1-\delta^2/2$, so that
\[
n\langle a\xi|\xi\rangle=\sum_{h\in F}\re\langle \pi_X(h)\xi|\xi\rangle \leq n-1+1-\delta^2/2=n-\delta^2/2,
\]
hence
$\langle a\xi|\xi\rangle \leq 1-\frac{\delta^2}{2n}$.
We thus have for every unit vector $\xi\in \ell^2(X)$ 
\[
\Vert \xi-a\xi\Vert\geq |\langle \xi-a\xi|\xi\rangle|\geq 1-\langle a\xi|\xi\rangle \geq \frac{\delta^2}{2n},
\]
and this proves that $1$ does not belong to the spectrum of $a$.
\hfill $\square$

\begin{remark}\label{rem1.3}
Let $G\curvearrowright X$ be an action. Then, as is well known, all orbits are infinite if and only if, for all non-empty finite sets $F_1,F_2\subset X$, there exists $g\in G$ such that $g\cdot F_1\cap F_2=\emptyset$. See for instance \citep[Lemma 4.4]{Kech}, \citep[Lemma 2.2]{KechTsan} or the nice proof of Lemma 2.4 in \cite{PoV}.

If it is not the case, the action $G\curvearrowright X$ is trivially amenable in the sense above: assume that there is an element $x_0\in X$ with finite orbit, i.e. $G\cdot x_0=\{g\cdot x:g\in G\}$ is a finite set. Then the functional $\f:\ell^\infty(X)\rightarrow\C$ defined by
\[
\f(a)=\frac{1}{|G\cdot x_0|}\sum_{y\in G\cdot x_0} a(y) \quad (a\in\ell^\infty(X))
\]
is a $G$-invariant state. Furthermore, the latter condition is equivalent to the fact that $\pi_X$ has a non-zero invariant vector.
By \citep[Lemma 4.2]{GlasMon}, this implies that, if $G$ has property (T) and if $G\curvearrowright X$ is amenable, then the latter is automatically trivially amenable.
\end{remark}
 
The next lemma is a slight extension of \citep[Lemma 2.1]{GlasMon}. The kind of positive $G$-maps that appear there play an important role in \cite{Anan}.

\begin{lemma}\label{lem1.4}
let $G$ be a group acting on sets $X$ and $Y$. If there exists a linear, unital, positive $G$-map $\Phi: \ell^\infty(Y)\rightarrow \ell^\infty(X)$ and if $G \curvearrowright X$ is amenable, then so is $G \curvearrowright Y$. It is the case whenever there is a $G$-map $\phi:X\rightarrow Y$.
\end{lemma}
\textsc{Proof.} If $\varphi$ is a $G$-invariant state on $\ell^\infty(X)$, then $\psi:=\varphi\circ\Phi$ is a $G$-invariant state on $\ell^\infty(Y)$.
\hfill $\square$

\par\vspace{2mm}
The following lemma is a consequence of \citep[Proposition 3.5]{Ros}.

\begin{lemma}\label{lem1.5}
Let $G$ be a group as above and let $G \curvearrowright X$ be an action of $G$ such that:
\begin{enumerate}
\item [(i)] every stabilizer $G _x=\{g\in G : g\cdot x=x\}$ is finite;
\item [(ii)] the action $G \curvearrowright X$ is amenable.
\end{enumerate}
Then $G$ is amenable.
\end{lemma}

\par\vspace{2mm}
We end the present section by recalling a few definitions and facts on von Neumann algebras. 
Let $M$ be a finite von Neumann algebra with separable predual, endowed with a finite, normal, faithful, normalized trace $\tau$, and let $1\in N\subset M$ be a von Neumann subalgebra of $M$. We denote by $\E_N$ the $\tau$-preserving conditional expectation from $M$ onto $N$, and by $x\mapsto \hat x$ the natural embedding of $M$ into $L^2(M,\tau)$ so that $\tau(x)=\langle x\hat 1|\hat 1\rangle$ for every $x\in M$. The adjoint map $x\mapsto x^*$ extends to $L^2(M,\tau)$  as an antilinear involution denoted by $J$; then the commutant $M'$ of $M$ in $B(L^2(M,\tau))$ is equal to $JMJ$.
Moreover, $\E_N$ extends to the orthogonal projection $e_N$ of $L^2(M,\tau)$ onto $L^2(N,\tau)$. 

Assume that $M$ is a factor. An automorphism $\theta$ of $M$ is \textit{outer} if, given $x\in M$, the condition:
\[
xy=\theta(y)x\quad\forall y\in M
\]
implies that $x=0$. This definition is equivalent to the usual one by \citep[Corollary 1.2]{Kall}, for instance. 

Even if it is always true in the case of finite factors, we assume throughout the article that all automorphisms are trace-preserving.

Let $G$ be a group and let $\sigma:G\rightarrow\mathrm{Aut}(M)$ be an action of $G$ on $M$. It is \textit{outer} if, for every $g\in G^*$, the automorphism $\sigma_g$ is outer. Recall also that $\sigma$ is \textit{ergodic} if the subalgebra of $\sigma$-invariant elements 
$M^\sigma:=\{x\in M: \sigma_g(x)=x\ \forall g\in G\}$ is equal to $\C$.

\par\vspace{2mm}

Let now $\omega$ be a free ultrafilter on $\N$ and let $M$ be a finite von Neumann algebra endowed with a finite trace as above. Then
\[
I_\omega=\{(a_n)\in\ell^\infty(\N,M): \lim_{n\to\omega}\Vert a_n\Vert_2=0\}
\]
is a closed two-sided ideal of the von Neumann algebra $\ell^\infty(\N,M)$ and the corresponding quotient algebra is denoted by $M^\omega$.  

We write $(a_n)^\omega=(a_n)+I_\omega$ for the equivalence class of $(a_n)$ in $M^\omega$, and we recall that $M$ embeds naturally into $M^\omega$, the image of $a\in M$ being the class of the constant sequence $(a,a,\ldots)$. The algebra $M^\omega$ is a finite von Neumann algebra and it is endowed with a natural faithful, normal trace $\tau^\omega$ given by
\[
\tau^\omega((x_n)^\omega)=\lim_{n\to\omega}\tau(x_n)\quad ((x_n)^\omega\in M^\omega).
\]

The relative commutant of $M$ in $M^\omega$ is sometimes denoted by $M_\omega$; every element of $M_\omega$ is represented by a bounded sequence $(x_n)$ such that $\lim_{n\to\omega}\Vert [x,x_n]\Vert_2=0$ for every $x\in M$, where, for all $a,b\in M$, $[a,b]:=ab-ba$. 

If $M$ is a type II$_1$ factor, then it is well known that $M^\om$ is also a type II$_1$ factor. 

\par\vspace{2mm}

Let $N$ be a von Neumann subalgebra of the finite von Neumann algebra $M$. For any free ultrafilter $\om$ on $\N$, it follows from \citep[Proposition 4.2.7]{GHJ} that 
the following diagram
\[
\begin{array}{ccl}
N'\cap M^\om & \subset & M^\om\\
\cup & & \cup\\
N'\cap N^\om & \subset & N^\om
\end{array}
\]
is a commuting square, \textit{i.e.} the restriction $\E$ of the conditional expectation $\E_{N^\om}$ to $N'\cap M^\om$ equals the trace-preserving conditional expectation onto $N'\cap N^\om$.

\par\vspace{2mm}

Any (trace-preserving) automorphism $\theta$ of $M$ extends naturally to an automorphism $\theta^\omega$ of $M^\omega$, whose restriction to $M_\omega$ is an automorphism of $M_\omega$. The automorphism $\theta$ is \textit{centrally trivial} if $\theta^\omega(x)=x$ for every $x\in M_\omega$. An action $\sigma$ of $G$ on $M$ 
is said to be \textit{centrally free} if $\sigma_g$ is not centrally trivial for every $g\not=1$. For all this, see for instance \citep{Connes} or \cite{FGL}. We denote by $(M^\om)^\sigma$ the fixed point subalgebra of $M^\om$ for the action $\sigma^\om$.

\par\vspace{2mm}

Finally, we denote by $L(G)$ the von Neumann algebra generated by $\lambda$, \textit{i.e.} it is the bicommutant of the group $\lambda(G)$ in $B(\ell^2(G))$. It is called the \textit{group von Neumann algebra} of $G$, and it is endowed with the natural trace $\tau$ defined by
\[
\tau(x)=\langle x\delta_1|\delta_1\rangle\quad (x\in L(G ))
\]
where $(\delta_g)_{g\in G}$ denotes the canonical basis of $\ell^2(G)$. Every element $x\in L(G )$ admits a Fourier series decomposition $x=\sum_{g\in G }x(g)\lambda(g)$ where $x(g)=\tau(x\lambda(g^{-1}))$ and $\sum_g |x(g)|^2=\Vert x\Vert_2^2$. If $H$ is a subgroup of $G$, then $L(H)$ identifies to the von Neumann subalgebra of $L(G)$ formed by all elements $y$ for which $y(g)=0$ for every $g\in G \smallsetminus H$. 

We simply write $\E_H$ (resp. $e_H$) instead of $\E_{L(H)}$ (resp. $e_{L(H)}$). 

Notice that $e_H$ is the multiplication operator by the characteristic function $\chi_H$ on $\ell^2(G)$.

By a slight abuse of notation, we'll still denote by $\tau$ the state on $B(\ell^2(G))$ given by
\[
\tau(x)=\langle x\delta_1|\delta_1\rangle\quad (x\in B(\ell^2(G)))
\]
even if it is not a tracial state on $B(\ell^2(G))$. 

\begin{remark}\label{rem2.5}
It is folklore that condition $(\star)$ of Remark \ref{rem2.3} is equivalent to the irreducibility of $L(H)$ in $L(G)$. Also, $L(G)$ is a factor if and only if $G$ is an icc group.
\end{remark}

\section{Relative inner amenability}\label{InnerAm}

We start by giving some equivalent conditions to relative inner amenability. They are inspired by the pioneering articles on inner amenability: \citep{BH}, \citep{Effros} and \cite{Paschke}. 

Let $H$ be a proper subgroup of a group $G$.
We denote by $C^*(\alpha(H))$ the C$^*$-algebra generated by $\alpha(H)$ in $B(\ell^2(G ))$. We observe moreover that:
\begin{itemize}
\item $\tau$ is a tracial state on $C^*(\alpha(H))$ since $\tau(\alpha(g))=1$ for every $g\in H$;
\item the projection $e_H$ commutes with all elements of $C^*(\alpha(H))$.
\end{itemize}

\begin{theorem}\label{thm2.1}
Let $H$ be a subgroup of a group $G$. Then the following conditions are equivalent:
\begin{enumerate}
\item [(1)] There exists an $H$-invariant state on $\ell^\infty(G \smallsetminus H)$ with respect to the conjugation action, i.e. $H$ is inner amenable relative to $G$.
\item [(2)] There exists a net $(\eta_n)$ of unit vectors in $\ell^1(G )^+$ such that $\supp(\eta_n)\subset G \smallsetminus H$ for every $n$, and
\[
\lim_{n}\Vert\alpha(h)\eta_n-\eta_n\Vert_1=0
\]
for every $h\in H$.
\item [(3)] There exists a net $(\xi_n)$ of unit vectors in $\ell^2(G )$ such that $\supp(\xi_n)\subset G \smallsetminus H$ for every $n$, and
\[
\lim_{n}\Vert\alpha(h)\xi_n-\xi_n\Vert_2=0
\]
for every $h\in H$, i.e. the trivial representation $1_H$ is weakly contained in $\pi_{G \smallsetminus H}$.
\item [(4)] There exists a net $(F_n)$ of non empty finite sets of $G \smallsetminus H$ such that
\[
\lim_{n}\frac{|F_n\bigtriangleup hF_nh^{-1}|}{|F_n|}=0
\]
for every $h\in H$.
\item [(5)] For every finite set $F\subset H$ and for every function $f:F\rightarrow [0,1]$ such that $\sum_{g\in F}f(g)=1$, the number $1$ belongs to the spectrum of $\alpha(f)(1-e_H)$.
\item [(6)] There exists a state $\f$ on $B(\ell^2(G ))$ such that $\f|_{C^*(\alpha(H))}=\tau|_{C^*(\alpha(H))}$ and $\f(e_H)=0$. In particular, $e_H\notin C^*(\alpha(H))$.
\end{enumerate}
\end{theorem}
\textsc{Proof.} Conditions (1) to (5) are equivalent by Lemma \ref{lem1.2}. The rest of our proof is inspired by \citep{Paschke}.

Let us prove that condition (3) implies condition (6): Let $(\xi_n)$ be a net as in (3), and let $\f_n$ be the corresponding state on $B(\ell^2(G ))$:
\[
\f_n(x)=\langle x\xi_n|\xi_n\rangle\quad\forall x\in B(\ell^2(G )).
\]
Then $\f_n(e_H)=0$ for every $n$ and $\f_n(\alpha(h))\to 1$ for every $h\in H$. Indeed, we have:
\[
\Vert\alpha(h)\xi_n-\xi_n\Vert_2^2=2(1-\re\langle\alpha(h)\xi_n|\xi_n\rangle)\to 0,
\]
which implies that $\langle\alpha(h)\xi_n|\xi_n\rangle\to 1$ since the sequence has modulus at most $1$.

Let $\f$ be an accumulation point of the net $(\f_n)$.
Then 
$\f(e_H)=0$ and $\f(a)=\tau(a)$ for every element $a$ in the $*$-algebra generated by $\alpha(H)$ since $\f(\alpha(h))=1=\tau(\alpha(h))$ for every $h\in H$.

Conversely, let us prove that condition (6) implies (1). If $\f$ is a state as in condition (6), one has $\f(\alpha(h))=1$ for every $h\in H$, hence
$\f(\alpha(h)x)=\f(x)=\f(x\alpha(h))$ for all $x\in B(\ell^2(G))$.
As in \citep{Paschke}, its restriction to $\ell^\infty(G\smallsetminus H)$ is a $H$-invariant state with respect to the action by conjugation.
\hfill $\square$

\begin{remark}\label{rem2.2}
Let $H$ be a subgroup of a group $G$ such that the projection $e_H$ does not belong to $C^*(\alpha(H))$. We do not know if this implies that $H$ is inner amenable relative to $G$.
\end{remark}

\begin{example}\label{ex2.5}
Let $H$ be an inner amenable group, set $G =H\times H$, embed $H$ into $G$ diagonally and let $\varphi:H^*\rightarrow G \smallsetminus H$ be defined by $\varphi(h)=(h,1)$ for $h\in H^*$. Then it is an $H$-map, hence the action on $G \smallsetminus H$ is amenable by Lemma \ref{lem1.4}.

Conversely, if $H$ is inner amenable relative to $H\times H$ with respect to the above embedding, then $H$ is inner amenable: indeed, let $\phi:H\times H\rightarrow H$ be defined by $\phi(g,h)=gh^{-1}$. Then it is an $H$-map, it maps $(H\times H)\smallsetminus H$ onto $H^*$,
and by Lemma \ref{lem1.4} again, the action $H\curvearrowright H^*$ by conjugation is amenable. We are grateful to Alain Valette for this observation.

Notice furthermore that the pair $H\subset H\times H$ satisfies condition $(\star)$ of Remark \ref{rem2.3} if and only if $H$ is an icc group.

We also observe that Y. Stalder presents in \cite{Stal} examples of $HNN$-extensions which are inner amenable and icc. These families contain non-amenable Baumslag-Solitar groups contrary to what was stated in \cite{BH}.
\hfill $\square$
\end{example}

\par\vspace{2mm}
Here are two easy hereditary properties of relative inner amenability.

\begin{example}\label{ex2.6}
(1) Suppose that $K<H<G$ and that $H$ is inner amenable relative to $G$. Then $K$ is inner amenable relative to $G$, too. Indeed, the natural inclusion of $G \smallsetminus H$ into $G \smallsetminus K$ is a $K$-map, and Lemma \ref{lem1.4} applies.\\
(2) Let $H_j$ be a subgroup of some group $G_j$, $j=1,2$. 
If $H_1$ is inner amenable relative to $G_1$ then the direct product group $H_1\times H_2$ is inner amenable relative to $G_1\times G_2$: apply condition (4) in Theorem \ref{thm2.1} for instance. Furthermore, the pair $H_1\times H_2\subset G_1\times G_2$ satisfies condition $(\star)$ in Remark \ref{rem2.3} if both pairs $H_1\subset G_1$ and $H_2\subset G_2$ do.
\hfill $\square$
\end{example}

The next example deals with semidirect products and will be needed in \ref{prop2.8}.

\begin{example}\label{ex2.7}
Let $A$ be an infinite group and let $\sigma:H\rightarrow \mathrm{Aut}(A)$ be an action of $H$ on $A$ which is amenable in the sense that there exists a $\sigma$-invariant mean on $A^*$. Then $H$ is inner amenable relative to the crossed product group $G =A\rtimes_\sigma H$. Indeed, realize $G$ as the direct product set $A\times H$ endowed with the composition law 
\[
(a_1,g_1)(a_2,g_2)=(a_1\sigma_{g_1}(a_2),g_1g_2).
\]
If $(\Phi_n)$ is a F\o lner net for the action $\sigma$ on $A^*$, set $F_n=\Phi_n\times\{1\}$ for every $n$. Then
$(F_n)$ is a F\o lner net for the action of $H$ by conjugation.

 Notice moreover that condition $(\star)$ in Remark \ref{rem2.3} is equivalent to the fact that the action $H\curvearrowright A^*$ has only infinite orbits, and this is equivalent to the fact that $H$ is not trivially inner amenable relative to $A\rtimes H$.
\hfill $\square$
\end{example}

\par\vspace{2mm}

We prove now that the family of restricted wreath products provides many examples of non-trivial relative inner amenable pairs of groups. This answers partly the question raised in the first section.

Assume that $H$ is a group that acts on some set $X$ and let $Z$ be any non-trivial group. Set 
\[
Z^{(X)}:=\{a:X\rightarrow Z: \supp(a)\ \textrm{is\ finite}\}.
\]
Then $H$ acts on $Z^{(X)}$ by Bernoulli shift: $\sigma_g(a)(x)=a(g^{-1}x)$. 

\begin{theorem}\label{prop2.8}
Let $H\curvearrowright X$ and $Z$ be as above. Then the action $H\curvearrowright Z^{(X)*}$ is amenable if and only if $H\curvearrowright X$ is amenable. If it is the case, $H$ is inner amenable relative to the restricted wreath product $G=Z\wr_XH$. 

Suppose moreover that $X$ is infinite. Then condition $(\star)$ of Remark \ref{rem2.3} holds if and only if all orbits of $H$ on $X$ are infinite. 
\end{theorem}
\textsc{Proof.}
Assume first that $H\curvearrowright X$ is amenable.
Choose some $z_0\in Z$, $z_0\not=1$. For every $x\in X$ let $\f_x\in Z^{(X)}$ be defined by
\[
\f_x(y)=\left\{
\begin{array}{ll}
1, & y\not=x\\
z_0, & y=x.
\end{array}\right.
\]
Since $\sigma_g(\f_x)=\f_{gx}$, as it is easily verified, the map $\f: x\mapsto \f_x$ is one-to-one and $H$-equivariant, and the claim follows from Lemma \ref{lem1.4}.

Assume now that $H\curvearrowright Z^{(X)*}$ is amenable. We claim first that $Z$ can be chosen to be equal to $\Z/2\Z$. Indeed, the map
$\psi:Z^{(X)*}\rightarrow (\Z/2\Z)^{(X)*}$ defined by
\[
\psi(a)(x)=\left\{
\begin{array}{cc}
1, & a(x)\not=1_Z\\
0, & a(x)=1_Z
\end{array}\right.
\]
is an $H$-map, hence the action $H\curvearrowright (\Z/2\Z)^{(X)*}$ is amenable by Lemma \ref{lem1.4}.

Thus let us assume that $Z=\Z/2\Z$. Then the corresponding action is closely related to the action by Bernoulli shift on the probability space $(Z^X,\mu^X)$ where $\mu^X$ denotes the product measure on $Z^X:=\prod_{x\in X}Z$ of the uniform probability measure on the finite set $Z$. 
It follows from \citep[Section 3]{KechTsan} that the permutation representation associated to the action $H\curvearrowright Z^{(X)*}$ 
is equivalent to 
the Koopman representation $\kappa_0$ of $H$ on $L^2(Z^X,\mu^X)_0:=\{\xi\in L^2(Z^X) : \int \xi d\mu^X=0\}$. By \citep[Theorem 1.2]{KechTsan}, the trivial representation of $H$ is weakly contained in $\pi_X$, and Lemma \ref{lem1.2} applies.

Assume now that $X$ is infinite. If $H\curvearrowright X$ has only infinite orbits, let us show that for any finite sets $F_1,F_2\subset Z^{(X)}\smallsetminus\{1\}$ there exists an element $g\in H$ such that $\sigma_g(F_1)\cap F_2=\emptyset$ (see Remark \ref{rem1.3}): let $S\subset X$ be a finite subset such that $\supp(a)\subset S$ for every $a\in F_1\cup F_2$, \textit{i.e.}, one has $a(x)=1$ for all $a\in F_1\cup F_2$ and $x\notin S$. Then there exists $g\in H$ such that $g^{-1}S\cap S=\emptyset$. If $a\in F_1$, there exists $x\in S$ such that $a(x)\not=1$. Then $\sigma_g(a)(gx)=a(x)\not=1$, so that if $\sigma_g(a)\in F_2$, then $x$ would belong to $g^{-1}S\cap S$, which is a contradiction. 

Conversely, if $H\curvearrowright X$ has some finite orbit $H\cdot x_0$, choose $z_0\in Z^*$ and set \[
a_0(x)=\left\{
\begin{array}{ll}
z_0 & x\in H\cdot x_0\\
1 & x\notin H\cdot x_0.
\end{array}\right.
\]
Then $\sigma_g(a_0)=a_0$ for every $g\in H$ and $a_0\not=1$.
\hfill $\square$

\begin{corollary}
For every group $H$ in the class $\mathcal B$ defined in Remark \ref{rem2.9}, there exists a group $G$ containing $H$ such that the latter is non-trivially inner amenable relative to $G$.
\end{corollary}

In some cases, relative inner amenability of $H$ in $G$ implies amenability of $H$; the following result is a consequence of Lemma \ref{lem1.5}.

\begin{proposition}\label{prop2.10}
Let $H\subset G$ be a pair of groups such that:
\begin{enumerate}
\item [(i)] $H$ is inner amenable relative to $G$;
\item [(ii)] for every $g\in G \smallsetminus H$, the subgroup $gHg^{-1}\cap H$ is finite, i.e. $H$ is an almost malnormal subgroup of $G$.
\end{enumerate}
Then $H$ is amenable.
\end{proposition}
\textsc{Proof.} The second condition implies that all stabilizers of the action $H\curvearrowright G \smallsetminus H$ are finite. Then the amenability of $H$ follows from Lemma \ref{lem1.5}.
\hfill $\square$

\par\vspace{2mm}
We end the present section with families of pairs $H\subset G$ that do not satisfy the relative inner amenability condition, even when $H$ has finite index in $G$.

\begin{example}\label{ex2.11}
Let $H$ be a group, and let $K$ be a non trivial group. Let $G =H*K$ be their free product. Let $C\subset G \smallsetminus H$ be the set of all elements $w=k_1g_1\ldots k_ng_n$ where $n\geq 1$, $k_i\in  K^*$ and $g_i\in H$ for every $i$, $w$ being in reduced form, so that $g_1,\ldots,g_{n-1}\in H^*$ if $n>1$. 
Then it is easily verified that, for $h,h'\in H$ and $w,w'\in C$ such that $hwh^{-1}=h'w'h'^{-1}$, one has $h=h'$ and $w=w'$. Furthermore, for every $g\in G\smallsetminus H$, there exists $w\in C$ and $h\in H$ such that $g=hwh^{-1}$. Thus
\[
G\smallsetminus H=\bigsqcup_{w\in C}\{hwh^{-1}: h\in H\}
\]
and the action of $H$ on each orbit is simply transitive, \textit{i.e.} it is equivalent to the action on $H$ by left translation.
Hence $H\curvearrowright G \smallsetminus H$ is amenable if and only if $H$ is amenable.
\hfill $\square$
\end{example}

\begin{example}\label{ex2.12}
The present example relies again on \cite{GlasMon}.
Let $H$ be an infinite group which has property (F) as discussed in Remark \ref{rem2.9}, but which does not have property (T). Such groups exist since property (F) is preserved under finite free products, by \cite[Section 4]{GlasMon}. Then $H$ is not inner amenable relative to any group $G \supset H$ for which the pair satisfies condition $(\star)$ of Remark \ref{rem2.3}. This shows that property (T) is not the only obstruction to the lack of non-trivial relative inner amenability.
\hfill $\square$
\end{example}

The next example shows the existence of pairs $H\subset G$ with $[G:H]<\infty$ such that $H$ is not inner amenable relative to $G$. It is based on a variant of the paradoxical decomposition of the free group that was used in \cite{Effros} in order to prove the non-inner amenability of the free group.

\begin{example}\label{ex5.1}
Let $H=\mathbb{F}_2$ be the non-abelian free group on two generators $a$ and $b$, and let $\theta$ be the order two automorphism of $H$ that exchanges $a$ and $b$. Denote by $G=H\rtimes\{1,\theta\}$ the corresponding semidirect product group. Then the action $H\curvearrowright G\smallsetminus H=H\times\{\theta\}$ is equivalent to the following action of $H$ on itself: $h\cdot x:=hx\theta(h^{-1})$ for all $h,x\in H$. In particular, for every $h\in H$, one has $a\cdot h=ahb^{-1}$ and $b^{\pm 1}\cdot h=b^{\pm 1}ha^{\mp 1}$.

Let then $E$ be the set of all elements of $H^*$ whose reduced form ends with an element of the form $b^m$ where $m$ is a non-zero integer. Then it is easily checked that: 
\begin{itemize}
\item $H=E\cup a\cdot E$;
\item the three sets $E$, $b\cdot E$ and $b^{-1}\cdot E$ are pairwise disjoint.
\end{itemize}
This implies that there cannot exist any invariant mean on $H$ with respect to the above action of $H$ and therefore that $H$ is not inner amenable relative to $G$.

Choosing any non-trivial inner amenable group $K$ and setting $H'=H\times K\subset G'=G\times K$, we get an example of a pair where $H'$ has finite index in $G'$, $H'$ is inner amenable, but it
is not inner amenable relative to $G'$. Indeed, as there is an obvious $H$-map from $G'\smallsetminus H'$ to $G\smallsetminus H$,
the relative inner amenability of $H'$ in $G'$ would imply that of $H$ in $G$ by Lemma \ref{lem1.4}.
\hfill $\square$
\end{example}

\begin{remark}
Let $H$ be an inner amenable group which is a subgroup of finite index of some group $G$. Example \ref{ex5.1} shows that $H$ is not necessarily inner amenable relative to $G$. However,
it is proved in \citep[Ajout]{BH} that $G$ itself is inner amenable. The proof consists in defining a $G$-invariant state $\f$ on $\ell^\infty(G^*)$ from an $H$-invariant mean $\mu$ on $H^*$. Let us remind its construction: we set, for every $f\in\ell^\infty(G^*)$,
\[
\f(f)=\frac{1}{[G:H]}\sum_{x\in G/H}\tilde{f}(x)
\]
where
\[
\tilde{f}(gH)=\int_H f(ghg^{-1})d\mu(h)\quad (gH\in G/H).
\]
The restriction $\psi$ of $\f$ to $\ell^\infty(G\smallsetminus H)$ is obviously $H$-invariant, but $\psi$ may be equal to zero. It is the case for instance when $H$ is normal in $G$.
\end{remark}

\section{Relative property gamma and non-Kazhdan groups}\label{propgamma}

As in the case of pairs of groups, it is straightforward to define a notion of relative property gamma: if $N\subset M$ is a pair of finite von Neumann algebras with separable preduals endowed with a finite trace $\tau$ as in Section 1, recall from Definition \ref{def1.1} that $N$ has \textit{property gamma relative to} $M$ if there exists a bounded sequence $(x_n)_{n\geq 1}\subset M$ such that $\Vert x_n\Vert_2=1$ and $\E_N(x_n)=0$ for every $n$, and such that $\Vert [x_n,y]\Vert_2\to 0$ as $n\to\infty$ for every $y\in N$. Because of the separability condition, it is equivalent to require that
\[
N'\cap M^\om\supsetneq N_\om (=N'\cap N^\om).
\]
We insist on the fact that our relative property gamma is distinct from the one in \cite{Bi1} and \cite{Bi2} (which can be expressed as the non-triviality of the relative commutant $M'\cap N^\om$), and that neither of them implies the other in an obvious way.

Any element of the relative commutant $N'\cap M$ obviously belongs to $N'\cap M^\om$. Thus, in order to avoid trivial cases, we assume from now on that $M$ is a type II$_1$ factor with separable predual and that $N$ is an irreducible subfactor of $M$.

\begin{remark}\label{rem3.1}
(1) The algebra $N'\cap M^\om$ has already been studied, for instance in \cite{FGL}. It follows from \citep[Lemma 3.5]{FGL} that, if $N$ is an irreducible subfactor of the type II$_1$ factor $M$, then either $N'\cap M^\omega=\C$ or $N'\cap M^\omega$ is diffuse, \textit{i.e.} it has no atoms.\\
(2) If $H$ is a countable icc subgroup of a countable icc group $G$ satisfying condition $(\star)$ of Remark \ref{rem2.3} so that $L(H)$ is an irreducible subfactor of the type II$_1$ factor $L(G)$, and if $L(H)$ has property gamma relative to $L(G )$, then it is obvious that $H$ is inner amenable relative to $G$. Example \ref{ex3.11} below (which is based on S. Vaes' example) will show that, as in the case of single groups, relative inner amenability is strictly weaker than relative property gamma. 
\end{remark}

Here is a first class of examples of irreducible factors $N\subset M$ which have the relative property gamma. In fact, it is the starting point of the present article and it is at the core of the main result of this section.

\begin{proposition}\label{prop3.2}
Let $H$ be an icc countable group that acts on some type II$_1$ factor $Q$ and such that the corresponding action $\sigma$ has the following properties:
\begin{enumerate}
\item [(i)] $\sigma$ is ergodic;
\item [(ii)] the fixed point algebra $(Q^\omega)^\sigma$ is diffuse.
\end{enumerate}
Then $L(H)$ has property gamma relative to $Q\rtimes_\sigma H$. 
\end{proposition}
\textsc{Proof.}
Indeed, first of all, $L(H)$ is an irreducible subfactor of $Q\rtimes_\sigma H$ since $H$ is icc and $\sigma$ is ergodic.

Next, set $M=Q\rtimes_\sigma H$ and $N=L(H)$. We observe that, if $x=\sum_g x(g)\lambda(g)\in M$, then $\E_N(x)=\sum_g\tau(x(g))\lambda(g)$. For any $x=(x_n)^\om\in (Q^\omega)^\sigma$, one has $[x,y]=0$ for every $y\in L(H)$ considered as an element of the crossed product $Q^\omega\rtimes_{\sigma^\omega}H\supset L(H)$
since the sequence $\sup_n\Vert x_n\Vert<\infty$ and
\[
\lim_{n\to\om}\Vert \lambda(g)x_n-x_n\lambda(g)\Vert_2=\lim_{n\to\om}\Vert \sigma_g(x_n)-x_n\Vert_2=0
\quad (g\in H).
\]
The algebra $(Q^\omega)^\sigma$ being diffuse, take a projection $e\in(Q^\omega)^\sigma$ such that $\tau^\omega(e)=1/2$. By \citep[Proposition 1.1.3]{Connes2}, we choose a representative $(e_n)\subset Q$ such that each $e_n$ is a projection and $\tau(e_n)=1/2$ for every $n$. Set $u_n=2e_n-1$ for every $n$. Then each $u_n$ is a unitary operator
and the sequence $(u_n)$ satisfies 
\[
\lim_{n\to\omega}\Vert [x,u_n]\Vert_2=0
\]
for every $x\in L(H)$, and $\E_{N}(u_n)=\tau(u_n)=0$ for every $n$.
\hfill $\square$

\begin{proposition}\label{prop3.3}
Let $N$ be an irreducible subfactor of the type II$_1$ factor $M$ with separable predual. Assume that $N$ is a full factor, \textit{i.e.} $N_\om=\C$. Then $N$ has property gamma relative to $M$ if and only if, for all $y_1,\ldots,y_m\in N$ and $\varepsilon>0$, there exists $u\in U(M)$ such that $\E_N(u)=0$ and 
\[
\Vert [y_j,u]\Vert_2\leq\varepsilon\quad (j=1,\ldots,m).
\]
\end{proposition}
\textsc{Proof.} As $N$ has separable predual, it suffices to show the existence of a unitary element $u\in U(N'\cap M^\om)$ such that $\E_{N^\om}(u)=0$. But $\E_{N^\om}(N'\cap M^\om)=N_\om=\C$, hence the restriction of $\E_{N^\om}$ to $N'\cap M^\om$ is equal to $\tau^\om$. As $N'\cap M^\om$ is diffuse, we choose a projection $e\in N'\cap M^\om$ whose trace $\tau^\om(e)=1/2$ and we set $u=2e-1$ as in Proposition \ref{prop3.2} above.
\hfill $\square$

\par\vspace{2mm}
As in the case of groups, property (T) for $N$ in the sense of \cite{CJ} implies that $N$ cannot have property gamma relative to $M$ in which it is an irreducible subfactor:

\begin{proposition}\label{prop3.4}
Let $N$ be a type II$_1$ subfactor of the type II$_1$ factor $M$. Suppose that $N$ has property (T).
Then $N'\cap M^\om=(N'\cap M)^\om$. In particular, if $N$ is irreducible in $M$, then 
$N'\cap M^\om=\C$. 
\end{proposition}
\textsc{Proof.} Recall from \cite{CJ} that $N$ has property (T) if and only if one can find $y_1,\ldots,y_m\in N$ and positive constants $\varepsilon$ and $K$ such that, for every $0<\delta<\varepsilon$, for every $N$-bimodule $H$ containing a unit vector $\xi$ such that 
\[
\max_j \Vert y_j\xi-\xi y_j\Vert<\varepsilon
\]
there exists a unit vector $\eta\in H$ such that $x\eta=\eta x$ for every $x\in N$ and $\Vert \eta-\xi\Vert<K\delta$. Let then $v=(v_n)^\om\in N'\cap M^\om$ be a unitary element. By assumption, we have
\[
\lim_{n\to\om}\Vert y_jv_n-v_ny_j\Vert_2=0
\]
for every $j$, hence there exists a set $A\in \om$ such that for every 
$n\in A$ there exists a unit vector $\eta_n\in L^2(M)$ such that $\Vert \eta_n-\hat{v}_n\Vert\leq \frac{1}{2n}$ and $u\eta_n u^*=\eta_n$ for all $u\in U(N)$. Then
\[
\Vert uv_nu^*-v_n\Vert_2\leq \Vert u\hat{v}_nu^*-u\eta_n u^*\Vert+\Vert \eta_n-\hat{v}_n\Vert\leq 1/n
\]
for every $u\in U(N)$ and every $n\in A$.

For $n\in A$, let then $K_n$ be the $\Vert\cdot\Vert_2$-closed convex hull of $\{uv_nu^*: u\in U(N)\}$ in $M$ and let $x_n$ be its element of minimal $\Vert\cdot\Vert_2$-norm. By uniqueness, one has $ux_nu^*=x_n$ for every $u\in U(N)$, hence $x_n\in N'\cap M$ and $\Vert x_n-v_n\Vert_2\leq 1/n$ for every $n\in A$. This implies that $v=(x_n)^\om\in (N'\cap M)^\om$.
\hfill $\square$
\par\vspace{2mm}

The family of pairs $N\subset M$ with relative property gamma that we present now constitutes the main result of the present section. We consider
$N=L(G)$ where $G$ is any icc group which does not have property (T). 

\par\vspace{2mm}
First of all, we need an auxiliary result which should be known to the specialists, but we did not find any appropriate proof in the literature.

\begin{lemma}\label{lem3.5}
Let $G$ be a group. Then the following conditions are equivalent:
\begin{enumerate}
\item [(1)] $G$ does not have Kazhdan's property (T);
\item [(2)] $G$ has a unitary representation $(\kappa,\Hh_\kappa)$ acting on a Hilbert space with the following properties:
\begin{enumerate}
\item [(a)] the trivial representation is weakly contained in $\kappa$;
\item [(b)] the representation $\kappa$ is \textbf{weakly mixing} in the sense of \cite{BR}: for every finite set $F\subset \Hh_\kappa$ and for every $\varepsilon>0$ there exists $g\in G$ such that 
\[
|\langle\kappa(g)\xi|\eta\rangle|<\varepsilon\quad (\xi,\eta\in F), 
\]
equivalently, $\kappa$ contains no non-trivial finite-dimensional subrepresentation \citep[Theorem 1.9]{BR};
\item [(c)] for every $g\in G^*$, there exists a unit vector $\xi\in \Hh_\kappa$ such that 
\[
|\langle\kappa(g)\xi|\xi\rangle|<1.
\]
\end{enumerate}
\end{enumerate}
Furthermore, if $G$ is countable, then $\Hh_\kappa$ can be chosen separable.
\end{lemma}
\textsc{Proof.} It is clear that if $G$ satisfies condition (2) then it cannot have property (T). Thus it remains to prove that (1) implies (2). 
As $G$ does not have property (T), there exists an unbounded, conditionally negative definite function $\psi:G\rightarrow \R_+$. Replacing it by $\psi+1-\delta_1$ if necessary, we assume that furthermore $\psi(g)\geq 1$ for every $g\in G^*$. For every real number $t>0$, let $(\pi_t,\Hh_t,\xi_t)$ be the cyclic representation associated to the positive definite function $\f_t=e^{-t\psi}$. Put then 
\[
\kappa=\bigoplus_{n\geq 1}\pi_{1/n}.
\]
It is easy to check that it satisfies conditions (a) and (c) above. As the direct sum of weakly mixing representations is obviously weakly mixing by \citep[Theorem 1.9]{BR}, it suffices to prove that every representation $\pi_t$ is weakly mixing, but this follows for instance directly from \citep[Lemma 4.4]{JolT}.
\hfill $\square$

\begin{theorem}\label{thm3.6}
Let $G$ be a countable icc group. Then $G$ does not have property (T) if and only if
there exists an action $\sigma$ of $G$ on the hyperfinite II$_1$ factor $R$ such that:
\begin{enumerate}
\item [(1)] $\sigma$ is weakly mixing, i.e. for every finite set $F\subset R$ and $\varepsilon>0$, one can find $g\in G$ such that
\[
(wm)\quad |\tau(\sigma_g(a)b)-\tau(a)\tau(b)|<\varepsilon\quad (a,b\in F);
\]
\item [(2)] $\sigma$ is centrally free and the fixed point algebra $(R_\om)^\sigma$ is of type II$_1$.
\end{enumerate}
In particular, if $G$ does not have property (T), then it admits an action $\sigma$ on $R$ such that $L(G)$ is an irreducible subfactor of the crossed product $R\rtimes_\sigma G$ and $L(G)$ has property gamma relative to $R\rtimes_\sigma G$.
\end{theorem}

As centrally free, weakly mixing actions of $G$ on $R$ imply that the crossed product algebra $R\rtimes_\sigma G$ is a type II$_1$ factor, we obtain the following consequence of Theorem \ref{thm3.6} and Proposition \ref{prop3.4}.

\begin{corollary}\label{cor3.7}
Let $G$ be a countable icc group. Then it does not have property (T) if and only if the type II$_1$ factor $L(G)$ can be embedded as an irreducible subfactor of some type II$_1$ factor $M$ so that $L(G)$ has property gamma relative to $M$.
\end{corollary}

Notice that if $G$ has the Haagerup property, then \citep[Theorem 2.3.4]{ccjjv} shows that it admits an action $\sigma$ on $R$ such that:
\begin{itemize}
\item $\sigma$ is strongly mixing and centrally free;
\item the fixed-point algebra $(R_\omega)^\sigma$ of the centralizer of $R$ in $R^\om$ is of type II$_1$.
\end{itemize}

\par\vspace{2mm}
\textsc{Proof of Theorem \ref{thm3.6}.} If $G$ has property (T), then $L(G)$ cannot be embedded into a II$_1$-factor $M$ as an irreducible subfactor with relative property gamma, by Proposition \ref{prop3.4}.

Conversely, let us assume that $G$ does not have property (T) and let us prove that it admits an action as stated.
The idea is to use the same construction as in Theorems 2.3.2 and 2.3.4 of \cite{ccjjv}, so that we only sketch the construction. We fix a representation $(\kappa,\Hh_\kappa)$ having all properties as stated in Lemma \ref{lem3.5}, we assume that the scalar product on $\Hh_\kappa$ is antilinear in the first variable, 
we set $\Hh=\ell^2(\N)\otimes \Hh_\kappa$ and $\pi=1\otimes\kappa$, which is a separable representation of $G$ which has the same properties as $\kappa$.
We realize the hyperfinite II$_1$-factor $R$ (with separable predual) as the von Neumann algebra generated by the 
Gel'fand-Naimark-Segal construction of the pair $(\mathrm{CAR}(\Hh),\tau)$ where $A:=\mathrm{CAR}(\Hh)$ is the C$^*$-algebra generated by $\{a(\xi): \xi\in \Hh\}$, $a:\Hh\rightarrow A$ being a linear isometry, and the following \textbf{canonical anticommutation relations} hold:
\begin{enumerate}
\item [(1)] $a^*(\xi)a(\eta)+a(\eta)a^*(\xi)=\langle\xi|\eta\rangle$
\item [(2)] $a(\xi)a(\eta)+a(\eta)a(\xi)=0$
\end{enumerate}
for all $\xi,\eta\in \Hh$. It is well known that $A$ is a uniformly hyperfinite C$^*$-algebra, and that the representation $\pi$ induces an action $\sigma$ of $G$ on $A$ characterized by
\[
\sigma_g(a(\xi))=a(\pi(g)\xi) \quad (g\in G, \xi\in \Hh).
\]
Moreover, there is a unique tracial state $\tau$ on $A$ such that
\[
\tau(a^*(\xi_m)\ldots a^*(\xi_1)a(\eta_1)\cdots a(\eta_n))=2^{-n}\delta_{n,m}\det(\langle\xi_i|\eta_j\rangle)
\]
for all $\xi_1,\ldots,\xi_m,\eta_1,\ldots,\eta_n\in \Hh$.

Then the action $\sigma$ extends to $R$ because it preserves $\tau$. 

In order to prove that $\sigma$ is weakly mixing, it suffices to check condition $(wm)$ on finite sets of elements of the form $a^*(\xi_m)\ldots a^*(\xi_1)a(\eta_1)\cdots a(\eta_n)$, where $\xi_1,\ldots,\xi_m$ and $\eta_1,\ldots,\eta_n$ belong to some finite set $F\subset \Hh$, and $n,m\leq L$ for some integer $L$. But condition (2b) in the previous lemma implies the existence of a sequence $(g_k)\subset G$ such that
\[
\max_{\xi,\eta\in F}|\langle\pi(g_k)\xi|\eta\rangle|\to 0
\]
as $k\to\infty$. Then 
the proof of Part (1) of \citep[Theorem 2.3.2]{ccjjv} adapts to prove that 
\[
\lim_{k\to\infty}\tau(\sigma_{g_k}(a)b)=\tau(a)\tau(b)
\]
for every $a=a^*(\xi_m)\ldots a^*(\xi_1)a(\eta_1)\ldots a(\eta_n)$, $b=a^*(\zeta_r)\ldots a^*(\zeta_1)a(\om_1\ldots a(\om_s)$, where
$\xi_1,\ldots\xi_m,\eta_1,\ldots,\eta_n,\zeta_1,\ldots,\zeta_r,\om_1,\ldots \om_s\in F$ and $m,n,r,s\leq L$.

Next, to prove that $\sigma$ is centrally free, we argue exactly as in \citep[Theorem 2.3.4]{ccjjv}: fix an element $g\in G^*$; there exists a unit vector $\eta\in \Hh_\kappa$ such that $|\langle\kappa(g)\eta|\eta\rangle|<1$. For all $n\geq 1$, set then
\[
e_n=a^*(\delta_n\otimes\eta)a(\delta_n\otimes\eta),
\]
so that $e_n$ is a projection of $A$ with trace $1/2$. Then it is easy to see that $e:=(e_n)^\om$ is a projection of the centralizing algebra $R_\om$, and by the above formulas, we see that
\[
\Vert \sigma_g(e_n)-e_n\Vert_2^2=\frac{1}{2}-\frac{1}{2}|\langle\eta,\kappa(g)\eta\rangle|^2>0
\]
for every $n$, proving that $\sigma^\om_g((e_n)^\om)\not=(e_n)^\om$.

Finally, let $(\eta_n)_{n\geq 1}$ be a sequence of unit vectors in $\Hh_\kappa$ such that, for every $g\in G$, $\Vert\kappa(g)\eta_n-\eta_n\Vert\to 0$ as $n\to\infty$. Set
$\xi_n=\delta_n\otimes\eta_n$ and $\zeta_n=2^{-1/2}(\delta_n+\delta_{n+1})\otimes\eta_n$ and define $e_n=a^*(\xi_n)a(\xi_n)$ and $f_n=a^*(\zeta_n)a(\zeta_n)$, so that $e=(e_n)^\om$ and $f=(f_n)^\om$ both belong to $(R_\om)^\sigma$. Then
\[
\Vert [e_n,f_n]\Vert_2^2=\frac{1}{4}
\]
for all $n$. The proof of \citep[Theorem 2.2.1]{Connes2} shows that this suffices to see that $(R_\om)^\sigma$ is noncommutative, and of type II$_1$.
\hfill $\square$

\begin{remark}\label{rem3.8}
Let $H$ be a group which does not have property (T). One could ask whether there exists a group $G$ containing $H$ such that the pair $(H,G)$ satisfies condition $(\star)$ of Remark \ref{rem2.3} and that $H$ is inner amenable relative to $G$.
In fact, Example \ref{ex2.12} shows that it is impossible if $H$ has property (F).
\end{remark}

\par\vspace{2mm}

Next, we discuss the relationship between relative inner amenability of an icc group $H$ in a countable group $G$ and relative property gamma of $L(H)$ in $L(G)$. As observed in Remark \ref{rem3.1}, if $H\subset G$ is a pair of icc groups such that $L(H)$ has property gamma relative to $L(G)$, then $H$ is inner amenable relative to $G$. As will be seen below, the converse is false, by a slight modification of the construction in \cite{Vaes}. 

But let us make the following observations first. The second one is partly analogous to Example \ref{ex2.5}.

\begin{proposition}\label{prop3.9}
(1) Let $H\subset G$ be infinite countable groups such that $H$ is not inner amenable relative to $G$. Then
\[
L(H)'\cap L(G)^\om=L(H)_\omega.
\]
In particular, if $H$ is neither inner amenable nor inner amenable relative to $G$, then
\[
L(H)'\cap L(G)^\om=\C.
\]
In this case, $L(H)$ is an irreducible subfactor of the full factor $L(G)$.\\
(2) Let $H$ be a countable, icc group. Embed $H$ into $H\times H$ diagonally as in Example \ref{ex2.5}. Then $L(H)$ is an irreducible subfactor of $L(H\times H)=L(H)\overline{\otimes}L(H)$, and if it has property gamma then it has property gamma relative to $L(H\times H)$.
\end{proposition}
\textsc{Proof.} (1) Set $M=L(G)$ and $N=L(H)$ for short.
Observe first that the pair $H\subset G$ satisfies condition $(\star)$ of Remark \ref{rem2.3}, so that $N'\cap M=Z(N)$.
 
Let $x=(x_n)^\omega\in N'\cap M^\omega$. Then $\E_{N^\omega}(x)=(\E_N(x_n))^\omega\in N'\cap N^\om=N_\omega$. If $\lim_{n\to\omega}\Vert x_n-\E_N(x_n)\Vert_2$ was strictly positive, one could find a sequence $(\xi_n)$ as in condition (3) of Theorem \ref{thm2.1}, which is impossible. 
Thus $x\in N_\omega$.

If furthermore $H$ is not inner amenable, then $N_\omega=\C$. This implies that $N$ is a full factor, and the rest is obvious.

(2) Observe first that, if $x=\sum_{g,h\in H}x(g,h)\lambda(g,h)\in L(H\times H)$, then $\E_H(x)=\sum_{h\in H} x(h,h)\lambda(h,h)$. Assume that $L(H)$ has property gamma. There exists a sequence $(u_n)\subset L(H)$ of unitary elements such that $\tau(u_n)=0$ for every $n$ and $\Vert [y,u_n]\Vert_2\to 0$ as $n\to\infty$ for every $y\in L(H)$. Set $v_n=u_n\otimes 1=\sum_h u_n(h)\lambda(h)\otimes 1\in L(H\times H)$. Then it is a sequence of unitary elements and $\E_H(v_n)=0$ for every $n$. It is straighforward to check that $\Vert [y,v_n]\Vert_2\to 0$ as $n\to\infty$ for every $y\in L(H)$.
\hfill $\square$

\begin{remark}
Let $H$ be an icc countable group, and embed it diagonally into $H\times H$ as in Proposition \ref{prop3.9}. If $L(H)$ has property gamma relative to $L(H\times H)$ then $H$ is inner amenable relative to $H\times H$, hence it is inner amenable by Example \ref{ex2.5}, but we do not know whether $L(H)$ has property gamma.
\end{remark}

\begin{example}\label{ex3.10}
(1) Let $H$ and $K$ be groups as in Example \ref{ex2.11}.
Thus, if we assume furthermore that $H$ is not inner amenable, then the pair $H\subset G$ satisfies the conditions of Proposition \ref{prop3.9} and $L(H)'\cap L(H*K)^\omega =\C$.\\
(2) Let $H$ be an icc group with property (F). Then it is not inner amenable, and it follows from Example \ref{ex2.12} that, if $G$ is any group containing $H$ so that the pair $H\subset G$ satisfies condition $(\star)$ of Remark \ref{rem2.3}, then $L(H)'\cap L(G  )^\om=\C$.
\hfill $\square$
\end{example}

In \cite{Vaes}, S. Vaes constructed a countable, icc, inner amenable group $H$ and he proved that the associated factor $L(H)$ is full, answering an open question raised already by E. Effros \cite{Effros} when he introduced the notion of inner amenability. 


\begin{example}\label{ex3.11}
We show that there exists a pair of icc groups $\tilde H\subset G  $ such that $\tilde H$ is inner amenable relative to $G  $,  but $L(\tilde H)$ does not have property gamma relative to $L(G  )$.

Let us recall Vaes' construction: choose an increasing sequence of distincts prime numbers $(p_n)_{n\geq 0}$ and set
\[
H_n:=\left(\Z/p_n\Z\right)\sp 3,
\quad
K:=\bigoplus_{n\geq 0}H_n,\quad
\textrm{and}\quad
K_N:=\bigoplus_{n\geq N}H_n
\]
for every non-negative integer $N$. Let $\Lambda:=SL(3,\Z)$ act in a natural way on $H_n$ for every $n$, and hence diagonally on $K$ and on $K_N$. Define an increasing sequence of groups $L_0<L_1<\ldots$ inductively by $L_0=K\rtimes\Lambda$, and $L_N\hookrightarrow L_{N+1}=L_N*_{K_N}(K_N\times\Z)$ for every $N$ (where $K_N<K<L_0<L_N$), and finally, set 
\[
L=\varinjlim L_N.
\]
Then set $\tilde H=L\times\Lambda$ which embeds naturally into $G=L\times L$. We claim that $\tilde H$ is inner amenable relative to $G$. Indeed, as $L$ is inner amenable by \cite{Vaes}, we choose some F\o lner sequence $(F_n)$ relative to the action by conjugation of $L$ on $L^*$. Then, as $\Lambda$ has property (T), we can assume that for every $n$, $F_n$ is not contained in $\Lambda$. Hence $\Phi_n:=\{1\}\times F_n\subset G\smallsetminus \tilde H$ for every $n$ and it is a F\o lner sequence for the action by conjugation of $\tilde H$ on $G\smallsetminus \tilde H$.

Moreover, by \citep[Lemma 2]{Vaes}, for every $(g_1,g_2)\in G \smallsetminus (K\times K)$, the set 
\[
\{(h_1,h_2)(g_1,g_2)(h_1,h_2)^{-1}: (h_1,h_2)\in\Lambda\times \Lambda\}
\]
is infinite, hence $L(G)\cap L(\Lambda\times\Lambda)'= L(K\times K)\cap L(\Lambda\times \Lambda)'$. 

Therefore, if $(x_n)\subset L(G)$ is a bounded sequence such that $\Vert x_ny-yx_n\Vert_2\to 0$ for every $y\in L(H\times\Lambda)$, then, since $\Lambda\times \Lambda$ has property (T), one has 
\[
\lim_{n\to\infty}\Vert x_n-\E_{L(K\times K)\cap L(\Lambda\times \Lambda)'}(x_n)\Vert_2=0.
\]
It follows from \citep[Lemma 3]{Vaes} that $L(K\times K)\cap L(\Lambda\times \Lambda)'$ is an infinite tensor product of suitable 2-dimensional algebras, and this implies by \cite{Vaes} that 
\[
\lim_{n\to\infty}\Vert \E_{L(K\times K)\cap L(\Lambda\times \Lambda)'}(x_n)-\tau(x_n)\Vert_2=0.
\]
It follows that $\lim_{n\to\infty}\Vert x_n-\tau(x_n)\Vert_2=0$, so we cannot have both $\Vert x_n\Vert_2=1$ and $\E_{\tilde H}(x_n)=0$ for all $n$.
\hfill $\square$
\end{example}

We conclude the present article with some remarks on the case where $N$ has finite (probabilistic) index in $M$ in the sense of \cite{PiPo}, \textit{i.e.} when there exists a positive constant $c$ such that $\E_N(x^*x)\geq cx^*x$ for every $x\in M$. Then the \textit{index} of $\E_N$ is
\[
\mathrm{Ind}(\E_N):=(\max\{c\geq 0 : \E_N(x^*x)\geq cx^*x, x\in M\})^{-1}.
\]
We choose the latter definition of index for three reasons:
\begin{itemize}
\item It is easy to define.
\item It makes sense for arbitrary conditional expectations.
\item It coincides with Jones' index in the case of finite factors.
\end{itemize}

For all this, see \cite{PiPo}, and especially \citep[Proposition 2.1]{PiPo}. 

For future use, we recall the following well-known fact; we give a quick proof for the reader's convenience.

\begin{lemma}\label{lem4.13}
Let $1\in A\subset B$ be von Neumann algebras for which there exists a conditional expectation $\E:B\rightarrow A$ with finite index. If $A$ is finite-dimensional, then so is $B$.
\end{lemma}
\textsc{Proof.} Let $c>0$ be such that $\E(x^*x)\geq cx^*x$ for every $x\in B$. By \citep[Properties 1.1.2]{PopaCBMS} and \citep[Proposition 4.1]{JolHAP}, it follows that $\E$ is automatically normal, \textit{i.e.} $\sigma$-weakly continuous.

If $B$ was infinite-dimensional, it would contain a sequence of non-zero, pairwise orthogonal projections $(e_n)$. Hence $(e_n)$ would tend $\sigma$-weakly to $0$, as well as $(\E(e_n))$. But since $A$ is finite-dimensional, this would imply that $\Vert \E(e_n)\Vert\to 0$, which is impossible since $\Vert\E(e_n)\Vert\geq c$ for every $n$.
\hfill $\square$

\par\vspace{2mm}
Thus, let $N$ be an irreducible subfactor of the type II$_1$ factor $M$ (with separable predual). Recall from Section \ref{prerequisites} that for any free ultrafilter $\om$ on $\N$, the following diagram
\[
\begin{array}{ccl}
N'\cap M^\om & \subset & M^\om\\
\cup & & \cup\\
N'\cap N^\om & \subset & N^\om
\end{array}
\]
is a commuting square. 

Assume henceforth that $N$ has finite Jones' index in $M$. Then so does $N^\om$ in $M^\om$ and $[M^\om:N^\om]=[M:N]=:k$ by \citep[Proposition 1.10]{PiPo}. Hence
the commuting square property implies that the conditional expectation $\E$ of $N'\cap M^\om$ onto $N'\cap N^\om$ satisfies the following inequality:
\[
\E(x^*x)\geq k^{-1}x^*x\quad (x\in N'\cap M^\om).
\]

Then we have the following alternative: either $N$ is a full factor, \textit{i.e.} the relative commutant $N'\cap N^\om$ is equal to $\C$, or $N$ has property gamma, and $N'\cap N^\om$ is a diffuse von Neumann algebra (\textit{i.e.} it has no atoms).

If $N$ is a full factor, then $N'\cap N^\om=\C$ and $N'\cap M^\om$ is finite-dimensional by Lemma \ref{lem4.13}, thus it is trivial by \citep[Lemma 3.5]{FGL}. Hence $N$ does not have property gamma relative to $M$. 

If $N$ has property gamma, then so does $M$ by \citep[Proposition 1.11]{PiPo} and both algebras $N'\cap N^\om$ and $N'\cap M^\om$ are diffuse by \citep[Lemma 3.5]{FGL}. The inclusion can be strict, and thus $N$ can have property gamma relative to $M$, as the next example shows.

\begin{example}\label{5.3}
 Let $N$ be a type II$_1$ factor with separable predual endowed with an outer action $\sigma$ of some finite group $G$ such that every automorphism $\sigma_g$ is approximately inner. (Notice that $N$ must have property gamma; for instance, the hyperfinite II$_1$ factor $R$ satisfies these conditions.)
 
 Put $M=N\rtimes_\sigma G$. Every element $x\in M$ has a unique expression as 
 \[
 x=\sum_{g\in G}x(g)\lambda(g)
 \]
 with $x(g)\in N$ for every $g$, with $\E_N(x)=x(1)$, and such that $\lambda(g)y\lambda(g^{-1})=\sigma_g(y)$ for all $g\in G$, $y\in N$.
 
Then $N$ is an irreducible subfactor of $M$, and we claim that it has property gamma relative to $M$. Indeed, fix $g\in G^*$ and let $(u_n)\subset U(N)$ be a sequence of unitary operators such that $\Vert\sigma_g(y)-u^*_nyu_n\Vert_2\to 0$ for every $y\in N$. Put $x_n=u_n\lambda(g)\in M$. Then $x_n$ is a unitary element of $M$ and $\E_N(x_n)=0$ for every $n$. Furthermore, $\Vert [x_n,y]\Vert_2\to 0$ as $n\to\infty$ for every $y\in N$ because
\begin{eqnarray*}
\Vert [x_n,y]\Vert_2 & = &
\Vert u_n\lambda(g)y-yu_n\lambda(g)\Vert_2\\
& = &
\Vert u_n\lambda(g)y\lambda(g^{-1})-yu_n\Vert_2\\
&=&
\Vert u_n\sigma_g(y)-yu_n\Vert_2.
\end{eqnarray*}
Thus the element $x=(x_n)^\om$ is a unitary element of $N'\cap M^\om$, and $\E_{N^\om}(x)=0$.
\hfill $\square$
\end{example}

\par\vspace{2mm}

\bibliographystyle{plain}
\bibliography{ref1}

\end{document}